\documentclass[12pt,reqno]{amsart}
\usepackage{amsmath, amsfonts, amssymb, amsthm, hyperref}
\usepackage{bm}
\allowdisplaybreaks[4]
\def\bg{\bigg}
\def\({\bg(}
\def\){\bg)}

\def\1{{\bf 1}}

\def\pmod #1{\ ({\rm{mod}}\ #1)}

\def\<{\langle}
\def\>{\rangle}
\theoremstyle{plain}
\newtheorem{theorem}{Theorem}[section]
\newtheorem{lemma}{Lemma}

\theoremstyle{definition}

\theoremstyle{remark}

\makeatletter
\@namedef{subjclassname@2020}{%
  \textup{2020} Mathematics Subject Classification}
\makeatother
 \vspace{4mm}

\begin{document}
\hbox{Preprint}
\medskip

\title[Proof of a conjecture on permutations]
{Proof of a conjecture on permutations}
\author{Yue-Feng She}
\address {(Yue-Feng She) Department of Applied Mathematics, Nanjing Forestry University, Nanjing
210037, People’s Republic of China}
\email{she.math@njfu.edu.cn}

\author{Xin-Qi Luo}
\address {(Xin-Qi Luo, corresponding author) School of Sciences, Changzhou Institute of Technology, Changzhou, 213032, People’s Republic of China}
\email{luoxq@czust.edu.cn}

\keywords{permutation, subset sum, combinatoric.
\newline \indent Supported by the Tian Yuan Mathematical Foundation (Grant No. 12526613).}
\begin{abstract}
Given a positive integer $n$, define a function on the symmetric group $S_n$ by
\[
F(\tau)=\sum_{k=1}^{n}k^2\tau(k)^2.
\]
Motivated by a conjecture of Zhi-Wei Sun, we investigate the residue classes attained by
$F(\tau)$ modulo $2n+1$. We prove that for every integer $n>4$, the set
$\{F(\tau):\tau\in S_n\}$ contains a complete residue system modulo $2n+1$.
The proof is based on a family of involutions whose values are controlled by subset sums of squares.
\end{abstract}
\maketitle

\section{Introduction}
\setcounter{lemma}{0}
\setcounter{theorem}{0}
\setcounter{equation}{0}
\setcounter{conjecture}{0}
\setcounter{remark}{0}
\setcounter{corollary}{0}

In this paper, intervals always mean intervals of integers. As usual, for each positive integer $n$, we denote by $S_n$ the symmetric group on $[1,n]$. Questions concerning permutations satisfying arithmetic restrictions have a long history in
combinatorics and number theory. Such problems often exhibit a surprising interaction between
algebraic structures and additive phenomena. Classical examples include the prime-circle conjecture posed by Filz \cite{Filz}, later collected in Guy's well-known book\cite{Guy}, which states that there exists $\tau\in S_n$ for each positive even number $n$ such that $\tau(1)+\tau(2),\tau(2)+\tau(3),\dots,\tau(n)+\tau(1)$ are all primes. This conjecture remains open and readers may refer to \cite{CFG,SLZ} for recent progress. 

Recently, Sun \cite{Sun1,Sun2,Sun3} proposed numerous conjectures concerning
permutations and additive combinatorics. For example, Sun conjectured that for any integer $n>5$, there is a permutation $\tau\in S_n$ such that
$$
\sum_{k=1}^{n-1}\frac{1}{\tau(k)\tau(k+1)}=1.
$$
This conjecture is confirmed by Han \cite{Han}.

Another natural direction is to study the value distribution of arithmetic functions defined on
the symmetric group. Specially, given a polynomial expression involving $k$ and $\tau(k)$, one may ask
whether the corresponding values cover all residue classes modulo a given integer.
Questions of this type connect permutation theory with additive combinatorics and the theory
of subset sums.

We define a function $F$ on $S_n$ by
\[
F(\tau)=\sum_{k=1}^{n}k^2\tau(k)^2.
\]
The purpose of this paper is to establish the following theorem, which is first conjectured by Sun
 \cite[Conjecture 11.8]{Sun2}.

\begin{theorem}
For every integer $n>4$, the set
$
\{F(\tau):\tau\in S_n\}
$
contains a complete residue system modulo $2n+1$.
\end{theorem}

This paper is organized as follows. In Chapter 2, we will prove a lemma on long intervals contained in subset-sum sets of squares. In Chapter 3, we will prove the main theorem.

\section{Subset sums of squares}
\setcounter{lemma}{0}
\setcounter{theorem}{0}
\setcounter{equation}{0}
\setcounter{conjecture}{0}
\setcounter{remark}{0}
\setcounter{corollary}{0}

For each positive integer $n$, define the subset-sum set of squares
\[
A(n):=\left\{\sum_{k\in I}k^2:\ I\subseteq [1,n] \right\},
\]
and
\[
\Sigma_n:=\sum_{k=1}^{n}k^2=\frac{n(n+1)(2n+1)}6.
\]
We also define
\[
L(n):=\max\{ \,|I|:I\subseteq A(n),\ I \text{ is an interval}\,\}.
\]
That is to say, $L(n)$ is the length of the longest interval in $A(n)$.

The following lemma is the key additive ingredient.

\begin{lemma}\label{L1}
For every integer $n\ge10$,
$
[129,\Sigma_n-129]\subseteq A(n).
$
Consequently,
$
L(n)=\Sigma_n-257.
$
\end{lemma}

\begin{proof}
We proceed by induction on $n$.

For $n=10$, a direct computation shows that $[129,\Sigma_{10}-129]$ is contained in $A(10)$.

Assume that
$
[129,\Sigma_n-129]\subseteq A(n)
$
for some $n\ge10$. Since
$$
A(n+1)=A(n)\cup \bigl(A(n)+(n+1)^2\bigr),
$$
both intervals
$$
[129,\Sigma_n-129]
$$
and
$$
[129+(n+1)^2,\Sigma_n-129+(n+1)^2]
$$
are contained in $A(n+1)$.

Furthermore,
\[
\Sigma_n-129-\bigl(129+(n+1)^2\bigr)\ge6
\]
for every $n\ge10$.
Hence these two intervals overlap, implying
\[
[129,\Sigma_{n+1}-129]\subseteq A(n+1).
\]

This completes the induction. What left is to justify why this interval is the longest interval in $A(n)$. 

One may easily verify that $128$ does not belong to $A(10)$ and $A(11)$. Then 128 does not belong to $A(n)$ for $n\geq10$ since $n^2>128$ when $n\geq12$. Furthermore, it is obvious that $k$ belongs to $A(n)$ if and only if $\Sigma_n-k$ belongs to $A(n)$. Hence $128$ and $\Sigma_n-128$ do not belong to $A(n)$. The conclusion follows from the observation that 
$$A(n)\subseteq [0,\Sigma_n]$$
and
$$
\Sigma_n-257>128
$$
for $n\geq10$.
\end{proof}

\section{Proof of the main theorem}
\setcounter{lemma}{0}
\setcounter{theorem}{0}
\setcounter{equation}{0}
\setcounter{conjecture}{0}
\setcounter{remark}{0}
\setcounter{corollary}{0}

\noindent{\it Proof of the Main Theorem.} We divide our proof into three cases.

\textbf{Case 1. $n$ is odd and $n\ge21$.}

For each
$
i\in\left[1, (n-1)/2\right],
$
define the transposition
\[
\tau_i=(\,i,\ n+1-i\,).
\]
For a subset
$
I\subseteq \left[1, (n-1)/2\right],
$
let $\tau_I$ denote the product of the transpositions $\tau_i$ $(i\in I)$.
Since the transpositions are pairwise disjoint, a straightforward calculation gives
\begin{align*}
F(\tau_I)
=&\sum_{k=1}^{n}k^4
+\sum_{i\in I}
\Bigl(
2i^2(n+1-i)^2-i^4-(n+1-i)^4
\Bigr)\\
=&\sum_{k=1}^{n}k^4
-4(n+1)^2
\sum_{i\in I}
\left(\frac{n+1}{2}-i\right)^2.
\end{align*}
Note that the map
\[
i\mapsto \frac{n+1}{2}-i
\]
permutes $[1,(n-1)/2]$. Hence
\[
\{F(\tau_I)\bmod (2n+1)\}
=
\left\{
C-4(n+1)^2x \pmod{2n+1}:
x\in A\!\left(\frac{n-1}{2}\right)
\right\},
\]
where
\[
C=\sum_{k=1}^{n}k^4.
\]
Since
\[
\gcd(4(n+1)^2,2n+1)=1,
\]
multiplication by $4(n+1)^2$ is a bijection modulo $2n+1$.
By Lemma \ref{L1},
\[
L\!\left(\frac{n-1}{2}\right)\ge2n+1.
\]
Therefore the residue classes represented by
$A((n-1)/2)$ already cover every element of $\mathbb Z/(2n+1)\mathbb Z$.
Hence
\[
\{F(\tau):\tau\in S_n\}
\]
contains all residue classes modulo $2n+1$.

\medskip

\textbf{Case 2. $n$ is even and $n\ge22$.}

For
$
i\in[1,n/2-1],
$
define
\[
\tau_i=(i,n-i).
\]
For a subset
$
I\subseteq \left[1, n/2-1\right],
$
let $\tau_I$ denote the product of the transpositions $\tau_i$ $(i\in I)$.
Arguing exactly as above, we obtain
\[
F(\tau_I)
=
\sum_{k=1}^{n}k^4
-
4n^2
\sum_{i\in I}
\left(\frac n2-i\right)^2.
\]
Consequently,
\[
\{F(\tau_I)\bmod (2n+1)\}
=
\left\{
C-4n^2x \pmod{2n+1}:
x\in A\(\frac{n}{2}-1\)
\right\}.
\]
Since
\[
\gcd(4n^2,2n+1)=1
\]
and
\[
L\(\frac{n}{2}-1\)\ge2n+1,
\]
the conclusion follows exactly as in Case 1.

\medskip

\textbf{Case 3. $5\le n\le20$.}

The remaining finitely many cases were verified by direct computation.

\end{document}